\documentclass{article}%
\usepackage{amsmath}
\usepackage{amsfonts}
\usepackage{amssymb}
\usepackage{graphicx}%
\providecommand{\U}[1]{\protect\rule{.1in}{.1in}}
\newtheorem{theorem}{Theorem}

\newtheorem{corollary}[theorem]{Corollary}

\newtheorem{example}[theorem]{Example}

\newtheorem{lemma}[theorem]{Lemma}

\newtheorem{remark}[theorem]{Remark}

\begin{document}

\begin{center}
{\LARGE On the normal approximation for random fields via martingale methods}

\bigskip

Magda Peligrad and Na Zhang

\bigskip
\end{center}

Department of Mathematical Sciences, University of Cincinnati, PO Box 210025,
Cincinnati, Oh 45221-0025, USA. \texttt{ }

email: peligrm@ucmail.uc.edu

email: zhangn4@mail.uc.edu

\begin{center}
Abstract

\bigskip
\end{center}

We prove a central limit theorem for strictly stationary random fields under a
sharp projective condition. The assumption was introduced in the setting of
random sequences by Maxwell and Woodroofe. Our approach is based on new
results for triangular arrays of martingale differences, which have interest
in themselves. We provide as applications new results for linear random fields
and nonlinear random fields of Volterra-type.

\bigskip


\bigskip


\section{\textbf{Introduction }}

Martingale methods are very important for establishing limit theorems for
sequences of random variables. The theory of martingale approximation,
initiated by Gordin (1969), was perfected in many subsequent papers. A random
field consists of multi-indexed random variables $(X_{u})_{u\in Z^{d}}$. The
main difficulty when analyzing the asymptotic properties of random fields, is
the fact that the future and the past do not have a unique interpretation.
Nevertheless, it is still natural to try to exploit the richness of the
martingale techniques. The main problem consists of the construction of
meaningful filtrations. In order to overcome this difficulty mathematicians
either used the lexicographic order or introduced the notion of commuting
filtration. The lexicographic order appears in early papers, such as in
Rosenblatt (1972), who pioneered the field of martingale approximation in the
context of random fields. An important result was obtained by Dedecker (1998)
who pointed out an interesting projective criteria for random fields, also
based on the lexicographic order. The lexicographic order leads to normal
approximation under projective conditions with respect to rather large,
half-plane indexed sigma algebras. In order to reduce the size of the
filtration used in projective conditions, mathematicians introduced the
so-called commuting filtrations. The traditional way for constructing
commuting filtrations is to consider random fields which are functions of
independent random variables. We would like to mention several remarkable
recent contributions in this direction by Gordin (2009), El Machkouri et al.
(2013), Wang and Woodroofe (2013), Voln\'{y} and Wang (2014), and Cuny et al.
(2015), who provided interesting martingale approximations in the context of
random fields. It is remarkable that Voln\'{y} (2015) imposed the ergodicity
conditions to only one direction of the stationary random field. Other recent
results involve interesting mixing conditions such as in the recent paper by
Bradley and Tone (2017).

In this paper we obtain a central limit theorem for random fields, for the
situation when the variables satisfy a generalized Maxwell-Woodroofe
condition. This is an interesting projective condition which defines a class
of random variables satisfying the central limit theorem and its invariance
principle, even in its quenched form. This condition is in some sense minimal
for this type of behavior as shown in Peligrad and Utev (2005). Its importance
was pointed out, for example, in papers by Maxwell and Woodroofe (2000), who
obtained a central limit theorem (CLT); Peligrad and Utev (2005) obtained a
maximal inequality and the functional form of the CLT; Cuny and Merlev\`{e}de
(2014) obtained the quenched form of this invariance principle. The
Maxwell-Woodroofe condition for random fields was formulated in Wang and
Woodroofe (2013), who also pointed out a variance inequality in the context of
commuting filtrations.

Compared to the paper of Wang and Woodroofe (2013), our paper has double
scope. First, to provide a central limit theorem under generalized
Maxwell-Woodroofe condition that extends the original result of Maxwell and
Woodroofe (2000) to random fields. Second, to use more general random fields
than Bernoulli fields. Our results are relevant for analyzing some statistics
based on repeated independent samples from a stationary process.

The tools for proving these results will consist of new theorems for
triangular arrays of martingales differences which have interest in
themselves. We present applications of our result to linear random fields and
nonlinear random fields, which provide new limit theorems for these structures.

Our results could also be formulated in the language of dynamical systems,
leading to new results in this field.

\section{Results}

Everywhere in this paper we shall denote by$\ ||\cdot||$ the norm in $L^{2}.$
By $\Rightarrow$ we denote the convergence in distribution. In the sequel
$[x]$ denotes the integer part of $x.$ As usual, $a\wedge b$ stands for the
minimum of $a$ and $b.$

Maxwell and Woodroofe (2000) introduced the following condition for a
stationary processes $(X_{i})_{i\in Z}$, adapted to a stationary filtration
$(\mathcal{F}_{i})_{i\in Z}:$%
\begin{equation}
\sum\nolimits_{k\geq1}\frac{1}{k^{3/2}}||E(S_{k}|\mathcal{F}_{1}%
)||<\infty,\text{ }S_{k}=\sum\nolimits_{i=1}^{k}X_{i}, \label{MW1}%
\end{equation}
and proved a central limit theorem for $S_{n}/\sqrt{n}$. In this paper we
extend this result to random fields.

For the sake of clarity we shall explain first the extension to random fields
with double indexes and, at the end, we shall formulate the results for
general random fields.

We shall introduce a stationary random field adapted to a stationary
filtration. For constructing a flexible filtration it is customary to start
with a stationary random field $(\xi_{n,m})_{n,m\in Z}$ and to introduce
another stationary random field $(X_{n,m})_{n,m\in Z}$ defined by
\begin{equation}
X_{n,m}=f(\xi_{i,j},i\leq n,j\leq m), \label{defx}%
\end{equation}
where $f$ is a measurable function. Note that $X_{n,m}$ is adapted to the
filtration $\mathcal{F}_{n,m}=\sigma(\xi_{i,j},i\leq n,j\leq m).$ As a matter
of fact $X_{n,m}=T^{n}S^{m}(X_{0,0})$ where for all $u$ and $v,$
$Tf(....x_{-1,v},x_{0,v})=f(....x_{0,v},x_{1,v})$ ($T$ is the vertical shift)
and $Sf(....x_{u,-1},x_{u,0})=f(....x_{u,0},x_{u,1})$ ($S$ is the horizontal shift).

We raise the question of normal approximation for stationary random fields
under projection conditions with respect to the filtration $(\mathcal{F}%
_{n,m})_{n,m\in Z}$. In several previous results involving various types of
projective conditions, the methods take advantage of the existence of
commuting filtrations, i.e.
\[
E(E(X|\mathcal{F}_{a,b})|\mathcal{F}_{u,v})=E(X|\mathcal{F}_{a\wedge u,b\wedge
v}).
\]
This type of filtration is induced by an initial random field $(\xi
_{n,m})_{n,m\in Z}$ of independent random variables, or, more generally can be
induced by stationary random fields $(\xi_{n,m})_{n,m\in Z}$ where only the
columns are independent, i.e. $\bar{\eta}_{m}=(\xi_{n,m})_{n\in Z}$ are
independent. This model often appears in statistical applications when one
deals with repeated realizations of a stationary sequence. We prove this
property in Lemma \ref{commuting} in the Appendix.

It is interesting to point out that commuting filtrations can be described by
the equivalent formulation: for $a\geq u$ we have
\begin{equation}
E(E(X|\mathcal{F}_{a,b})|\mathcal{F}_{u,v})=E(X|\mathcal{F}_{u,b\wedge v}).
\label{pcf}%
\end{equation}
This follows from this Markovian-type property, see for instance Problem 34.11
in Billingsley (1995).

Our main result is the following theorem which is an extension of the CLT in
Maxwell and Woodroofe (2000) to random fields. Below we use the notation%
\[
S_{k,j}=\sum\nolimits_{u,v=1}^{k,j}X_{u,v}.
\]

\begin{theorem}
\label{indeprows0}Define $(X_{n,m})_{n,m\in Z}$ by (\ref{defx}) and assume
that (\ref{pcf}) holds. Assume that the following projective condition is
satisfied
\begin{equation}
\sum\nolimits_{j,k\geq1}\frac{1}{j^{3/2}k^{3/2}}||E(S_{j,k}|\mathcal{F}%
_{1,1})||<\infty. \label{MW}%
\end{equation}
In addition assume that the vertical shift $T$ is ergodic. Then there is a
constant $c$ such that%
\[
\frac{1}{n_{1}n_{2}}E(S_{n_{1},n_{2}}^{2})\rightarrow c^{2}\text{ as }%
\min(n_{1},n_{2})\rightarrow\infty
\]
and
\begin{equation}
\frac{1}{\sqrt{n_{1}n_{2}}}S_{n_{1},n_{2}}\Rightarrow N(0,c^{2})\text{ as
}\min(n_{1},n_{2})\rightarrow\infty. \label{cltmartf}%
\end{equation}

\end{theorem}

By simple calculations involving the properties of conditional expectation we
obtain the following corollary.

\begin{corollary}
\label{indeprows} Assume the following projective condition is satisfied
\begin{equation}
\sum\limits_{j,k\geq1}\frac{1}{j^{1/2}k^{1/2}}||E(X_{j,k}|\mathcal{F}%
_{1,1})||<\infty, \label{mixc1}%
\end{equation}
and $T$ is ergodic. Then there is a constant $c$ such that the CLT in
(\ref{cltmartf}) holds.
\end{corollary}

The results are easy to extend to general random fields $(X_{\mathbf{u}%
})_{\mathbf{u\in Z}^{d}}$ introduced in the following way. We start with a
stationary random field $(\xi_{\mathbf{n}})_{\mathbf{n\in}Z^{d}}$ and
introduce another stationary random field $(X_{\mathbf{n}})_{\mathbf{n\in
}Z^{d}}$ defined by $X_{\mathbf{k}}=f(\xi_{\mathbf{j}},\ \mathbf{j}%
\leq\mathbf{k}),$ where $f$ is a measurable function and $\mathbf{j}%
\leq\mathbf{k}$ denotes $j_{i}\leq k_{i}$ for all $i$. Note that
$X_{\mathbf{k}}$ is adapted to the filtration $\mathcal{F}_{\mathbf{k}}%
=\sigma(\xi_{\mathbf{u}},\mathbf{u}\leq\mathbf{k}).$ As a matter of fact
$Y_{\mathbf{k}}=T_{1}T_{2}...T_{d}(Y_{\mathbf{0}})$ where $T_{i}$ are the
shift operators.

In the next theorem we shall consider commuting filtrations in the sense that
for $a\geq u\in R^{1},$ $\mathbf{b},\mathbf{v}\in R^{d-1}$ we have
\[
E(E(X|\mathcal{F}_{a,\mathbf{b}})|\mathcal{F}_{u,\mathbf{v}})=E(X|\mathcal{F}%
_{u,\mathbf{b}\wedge\mathbf{v}}).
\]
For example, this kind of filtration is induced by stationary random fields
$(\xi_{n,\mathbf{m}})_{n\mathbf{\in}Z,\mathbf{m}\in Z^{d}}$ such that the
variables $\eta_{\mathbf{m}}=(\xi_{n,\mathbf{m}})_{n\mathbf{\in}Z^{1}}$ are
independent, $\mathbf{m}\in Z^{d-1}$. All the results extend in this context
via mathematical induction. Below, $|\mathbf{n}|=n_{1}\cdot...\cdot n_{d}.$

\begin{theorem}
\label{general}Assume that $(X_{\mathbf{u}})_{\mathbf{u\in Z}^{d}}$ and
$(\mathcal{F}_{\mathbf{u}})_{\mathbf{u\in Z}^{d}}$ are as above and assume
that the following projective condition is satisfied
\[
\sum\nolimits_{\mathbf{u}\geq\mathbf{1}}\frac{1}{|\mathbf{u}|^{3/2}%
}||E(S_{\mathbf{u}}|\mathcal{F}_{\mathbf{1}})||<\infty.
\]
In addition assume that $T_{1}$is ergodic. Then there is a constant $c$ such
that%
\[
\frac{1}{|\mathbf{n}|}E(S_{\mathbf{n}}^{2})\rightarrow c^{2}\text{ as }%
\min(n_{1},...,n_{d})\rightarrow\infty
\]
and
\begin{equation}
\frac{1}{\sqrt{|\mathbf{n}|}}S_{\mathbf{n}}\Rightarrow N(0,c^{2})\text{ as
}\min(n_{1},...,n_{d})\rightarrow\infty. \label{clt2}%
\end{equation}

\end{theorem}

\begin{corollary}
\label{indepWW}Assume that%
\begin{equation}
\sum\nolimits_{\mathbf{u}\geq\mathbf{1}}\frac{1}{|\mathbf{u}|^{1/2}%
}||E(X_{\mathbf{u}}|\mathcal{F}_{\mathbf{1}})||<\infty\label{mixg}%
\end{equation}
and $T_{1}$ is ergodic. Then the CLT in (\ref{clt2}) holds.
\end{corollary}

Corollary \ref{indepWW} above shows that Theorem 1.1 in Wang and Woodroofe
(2013)\ holds for functions of random fields which are not necessarily
functions of i.i.d.

\bigskip

We shall give examples providing new results for linear and Volterra random
fields. For simplicity, they are formulated in the context of functions of i.i.d.

\begin{example}
\label{exinear}(Linear field) Let $(\xi_{\mathbf{n}})_{\mathbf{n}\in Z^{d}}$
be a random field of independent, identically distributed random variables
which are centered and have finite second moment. Define%
\[
X_{\mathbf{k}}=\sum_{\mathbf{j}\geq\mathbf{0}}a_{\mathbf{j}}\xi_{\mathbf{k}%
-\mathbf{j}}.
\]
Assume that $\sum_{\mathbf{j}\geq\mathbf{0}}a_{\mathbf{j}}^{2}<\infty$ and%
\begin{equation}
\sum\limits_{\mathbf{j}\geq\mathbf{1}}\frac{|b_{\mathbf{j}}|}{|\mathbf{j}%
|^{3/2}}<\infty\text{ where }b_{\mathbf{j}}^{2}=\sum\limits_{\mathbf{i}%
\geq\mathbf{0}}(\sum\limits_{\mathbf{u}=\mathbf{1}}^{\mathbf{j}}%
a_{\mathbf{u}+\mathbf{i}})^{2}. \label{linear}%
\end{equation}
Then the CLT in (\ref{clt2}) holds.
\end{example}

Let us mention how this example differs from other results available in the
literature. Example 1 in El Machkouri et al. (2013) contains a CLT under the
condition $\sum_{\mathbf{u}\geq\mathbf{0}}|a_{\mathbf{u}}|<\infty.$ If we take
for instance for $u_{i}$ positive integers
\[
a_{u_{1},u_{2},...,u_{d}}=\prod\limits_{i=1}^{d}(-1)^{u_{i}}\frac{1}%
{\sqrt{u_{i}}\log u_{i}},
\]
then $\sum_{\mathbf{u}\in Z^{2}}|a_{\mathbf{u}}|=\infty.{\mathcal{\ }}%
$Furthermore, condition (\ref{mixg}), which was used in this context by Wang
and Woodroofe (2013), is not satisfied but condition (\ref{linear}) holds.

\bigskip

Another class of nonlinear random fields are the Volterra processes, which
plays an important role in the nonlinear system theory.

\begin{example}
\label{Volterra}(Volterra field) Let $(\xi_{\mathbf{n}})_{\mathbf{n}\in Z^{d}%
}$ be a random field of independent random variables identically distributed
centered and with finite second moment. Define%
\[
X_{\mathbf{k}}=\sum_{(\mathbf{u},\mathbf{v)}\geq(\mathbf{0},\mathbf{0}%
)}a_{\mathbf{u},\mathbf{v}}\xi_{\mathbf{k-u}}\xi_{\mathbf{k-v}},
\]
where $a_{\mathbf{u},\mathbf{v}}$ are real coefficients with $a_{\mathbf{u}%
,\mathbf{u}}=0$ and $\sum_{\mathbf{u,v}\geq\mathbf{0}}a_{\mathbf{u,v}}%
^{2}<\infty.$ Denote%
\[
c_{\mathbf{u},\mathbf{v}}(\mathbf{j})=\sum\limits_{\mathbf{k}=\mathbf{1}%
}^{\mathbf{j}}a_{\mathbf{k+u,k+v}}%
\]
and assume
\[
\sum\limits_{\mathbf{j}\geq1}\frac{|b_{\mathbf{j}}|}{|\mathbf{j}|^{3/2}%
}<\infty\text{ where }b_{\mathbf{j}}^{2}=\sum\limits_{\mathbf{u}%
\geq\mathbf{0,v}\geq\mathbf{0,u}\neq\mathbf{v}}(c_{\mathbf{u},\mathbf{v}}%
^{2}(\mathbf{j})+c_{\mathbf{u},\mathbf{v}}(\mathbf{j})c_{\mathbf{v}%
,\mathbf{u}}(\mathbf{j})).
\]
Then the CLT in (\ref{clt2}) holds.
\end{example}

\begin{remark}
In examples \ref{exinear} and \ref{Volterra} the fields are Bernoulli.
However, we can take as innovations the random field $(\xi_{n,m})_{n.m\in Z}$
having as columns independent copies of a stationary and ergodic martingale
differences sequence.
\end{remark}

\section{Proofs}

In this section we gather the proofs. They are based on a new result for a
random field consisting of triangular arrays of row-wise stationary martingale
differences, which allows us to find its asymptotic behavior by analyzing the
limiting distribution of its columns.

\begin{theorem}
\label{main}Assume that for each $n$ fixed $(D_{n,k})_{k\in Z}$ forms a
stationary martingale difference sequence adapted to the stationary nested
filtration $(\mathcal{F}_{n,k})_{k\in Z}$ and the family $(D_{n,1}^{2}%
)_{n\geq1}$ is uniformly integrable. In addition assume that for all $m\geq1$
fixed, $(D_{n,1},...,D_{n,m})_{n\geq1}$ converges in distribution to
$(L_{1},L_{2},...,L_{m}),$ and%
\begin{equation}
\frac{1}{m}\sum_{j=1}^{m}L_{j}^{2}\rightarrow c^{2}\text{ in }L^{1}\text{ as
}m\rightarrow\infty. \label{erglimit}%
\end{equation}
Then
\[
\frac{1}{\sqrt{n}}\sum\nolimits_{k=1}^{n}D_{n,k}\Rightarrow cZ\text{ as
}n\rightarrow\infty,
\]
where $Z$ is a standard normal variable.
\end{theorem}

\textbf{Proof of Theorem \ref{main}}. For the triangular array $(D_{n,k}%
/\sqrt{n})_{k\geq1},$ we shall verify the conditions of Theorem \ref{martCLT},
given for convenience in the Appendix. Note that for $\varepsilon>0$ we have%
\begin{equation}
\frac{1}{n}E(\max_{1\leq k\leq n}D_{n,k}^{2})\leq\varepsilon^{2}+E(D_{n,1}%
^{2}I(|D_{n,1}|>\varepsilon\sqrt{n})) \label{ineq0}%
\end{equation}
and, by the uniformly integrability of $(D_{n,1}^{2})_{n\geq1}$, we obtain:
\[
\lim_{n\rightarrow\infty}E(D_{n,1}^{2}I(|D_{n,1}|>\varepsilon\sqrt{n}))=0.
\]
Therefore, by passing to the limit in inequality (\ref{ineq0}),\ first with
$n\rightarrow\infty$ and then with $\varepsilon\rightarrow0,$ the first
condition of Theorem \ref{martCLT} is satisfied. The result will follow from
Theorem \ref{martCLT} if we can show that%
\[
\frac{1}{n}\sum_{j=1}^{n}D_{n,j}^{2}\rightarrow^{L^{1}}c^{2}\text{ }\ \text{as
}n\rightarrow\infty\text{.}%
\]

To prove it, we shall apply the following lemma to the sequence $(D_{n,k}%
^{2})_{k\in Z}$ after noticing that, under our assumptions, for all $m\geq1$
fixed, $(D_{n,1}^{2},...,D_{n,m}^{2})_{n\geq1}$ converges in distribution to
$(L_{1}^{2},L_{2}^{2},...,L_{m}^{2}).$

\begin{lemma}
\label{LLN}Assume that the triangular array of random variables $(X_{n,k}%
)_{k\in Z}$ is row-wise stationary and $(X_{n,1})_{n\geq1}$ is a uniformly
integrable family. For all $m\geq1$ fixed, $(X_{n,1},...,X_{n,m})_{n\geq1}$
converges in distribution to $(X_{1},X_{2},...,X_{m})$ and
\begin{equation}
\frac{1}{m}\sum_{u=1}^{m}X_{u}\rightarrow c\text{ in }L^{1}\text{ as
}m\rightarrow\infty. \label{conv1}%
\end{equation}
Then
\[
\frac{1}{n}\sum_{u=1}^{n}X_{n,u}\rightarrow c\text{ in }L^{1}\text{ as
}n\rightarrow\infty.
\]

\end{lemma}

\textbf{Proof of Lemma \ref{LLN}.} Let $m\geq1$ be a fixed integer and define
consecutive blocks of indexes of size $m$, $I_{j}(m)=\{(j-1)m+1,...,mj)\}.$ In
the set of integers from $1$ to $n$ we have $k_{n}=k_{n}(m)=[n/m]$ such blocks
of integers and a last one containing less than $m$ indexes. Practically, by
the stationarity of the rows and by the triangle inequality, we write%
\begin{gather}
\frac{1}{n}E|\sum\nolimits_{u=1}^{n}(X_{n,u}-c)|\leq\label{ineq}\\
\frac{1}{n}\sum\nolimits_{j=1}^{k_{n}}E|\sum\nolimits_{k\in I_{j}(m)}%
(X_{n,k}-c)|+\frac{1}{n}E|\sum\nolimits_{u=k_{n}m+1}^{n}(X_{n,u}%
-c)|\nonumber\\
\leq\frac{1}{m}E|\sum\nolimits_{u=1}^{m}(X_{n,u}-c)|+\frac{m}{n}%
E|X_{n,1}-c|.\nonumber
\end{gather}
Note that, by the uniform integrability of $(X_{n,1})_{n\geq1},$ we have
\[
\limsup_{n\rightarrow\infty}\frac{m}{n}E|X_{n,1}-c|\leq\limsup_{n\rightarrow
\infty}\frac{m}{n}\left(  E|X_{n,1}|+|c|\right)  =0.
\]
Now, by the continuous function theorem and by our conditions, for $m$ fixed,
we have the following convergence in distribution:
\[
\frac{1}{m}\sum\nolimits_{u=1}^{m}(X_{n,u}-c)\Rightarrow\frac{1}{m}%
\sum\nolimits_{u=1}^{m}(X_{u}-c).
\]
In addition, by the uniform integrability of $(X_{n,k})_{n}$ and by the
convergence of moments theorem associated to convergence in distribution, we
have\
\[
\lim_{n\rightarrow\infty}\frac{1}{m}E|\sum\nolimits_{u=1}^{m}(X_{n,u}%
-c)|=\frac{1}{m}E|\sum\nolimits_{u=1}^{m}(X_{u}-c)|,
\]
and by assumption (\ref{conv1}) we obtain%
\[
E|\frac{1}{m}\sum\nolimits_{u=1}^{m}(X_{u}-c)|\rightarrow\text{ }0\text{ as
}m\rightarrow\infty.
\]
The result follows by passing to the limit in (\ref{ineq}), letting first
$n\rightarrow\infty$ followed by $m\rightarrow\infty$. $\ \square$

\bigskip

When we have additional information about the type of the limiting
distribution for the columns the result simplifies.

\begin{corollary}
\label{cormain}If in Theorem \ref{main}\ the limiting vector $(L_{1}%
,L_{2},...,L_{m})$ is stationary Gaussian, then condition (\ref{erglimit})
holds and
\[
\frac{1}{\sqrt{n}}\sum\nolimits_{k=1}^{n}D_{n,k}\Rightarrow cZ\text{ as
}n\rightarrow\infty,
\]
where $Z$ is a standard normal variable and $c$ can be identified by%
\[
c^{2}=\lim_{n\rightarrow\infty}E(D_{n,1}^{2}).
\]

\end{corollary}

\textbf{Proof}. We shall verify the conditions of Theorem \ref{main}. Note
that, by the martingale property, we have\ that $\mathrm{cov}(D_{n,1}%
,D_{n,k})=0.\ $Next, by the condition of uniform integrability, by passing to
the limit we obtain $\mathrm{cov}(L_{1},L_{k})=0.$ Therefore, the sequence
$(L_{m})_{m}$ is an i.i.d. Gaussian sequence of random variables and condition
(\ref{erglimit}) holds. $\ \square$

\bigskip

In order to prove Theorem \ref{indeprows0} we start by pointing out an upper
bound for the variance given in Corollary 7.2 in Wang and Woodroofe (2013). It
should be noticed that to prove it, the assumption that the random field is
Bernoulli is not needed.

\begin{lemma}
\label{var0}Define $(X_{n,m})_{n,m\in Z}$ by (\ref{defx}) and assume that
(\ref{pcf}) holds. Then,\textbf{ }there is a universal constant $C$ such that
\[
\frac{1}{\sqrt{nm}}||S_{n,m}||\leq C\sum\nolimits_{i,j\geq1}\frac
{1}{(ji)^{3/2}}||E(S_{j,i}|\mathcal{F}_{1,1})||.
\]

\end{lemma}

By applying the triangle inequality, the contractivity property of the
conditional expectation and changing the order of summations we easily obtain
the following corollary.

\begin{corollary}
\label{var}Under the conditions of Lemma \ref{var0}\textbf{ }there is a
universal constant $C$ such that
\[
\frac{1}{\sqrt{nm}}||S_{n,m}||\leq C\sum\nolimits_{i,j\geq1}\frac
{1}{(ji)^{1/2}}||E(X_{j,i}|\mathcal{F}_{1,1})||\ .
\]

\end{corollary}

\textbf{ }

\textbf{Proof of Theorem \ref{indeprows0}}.

We shall develop the "small martingale method" in the context of random
fields. To construct a row-wise stationary martingale approximation we shall
introduce a parameter. Let $\ell$ be a fixed positive integer and denote
$k=[n_{2}/\ell]$. We start the proof by dividing the variables in each line in
blocks of size $\ell$ and making the sums in each block. Define
\[
X_{j,i}^{(\ell)}=\frac{1}{\ell^{1/2}}\sum_{u=(i-1)\ell+1}^{i\ell}%
X_{j,u}\;,i\geq1.
\]
Then, for each line $j$ we construct the stationary sequence of martingale
differences $(Y_{j,i}^{(\ell)})_{i\in Z}$ defined by
\[
Y_{j,i}^{(\ell)}=X_{j,i}^{(\ell)}-E(X_{j,i}^{(\ell)}|\mathcal{F}%
_{j,i-1}^{(\ell)}),
\]
where $\mathcal{F}_{j,k}^{(\ell)}=\mathcal{F}_{j,k\ell}$. Also, we consider
the triangular array of martingale differences $(D_{n_{1},i}^{(\ell)}%
)_{i\geq1}$ defined by%
\[
D_{n_{1},i}^{(\ell)}=\frac{1}{\sqrt{n_{1}}}\sum\nolimits_{j=1}^{n_{1}}%
Y_{j,i}^{(\ell)}.
\]
In order to find the limiting distribution of $(\sum\nolimits_{i=1}%
^{k}D_{n_{1},i}^{(\ell)}/\sqrt{k})_{k}$ when $\min(n_{1},,k)\rightarrow
\infty,$ we shall apply Corollary \ref{cormain}. It is enough to show that
\[
(D_{n_{1},1}^{(\ell)},...,D_{n_{1},N}^{(\ell)})\Rightarrow(L_{1},...,L_{N}),
\]
where $(L_{1},...,L_{N})$ is stationary Gaussian and $[(D_{n_{1},1}^{(\ell
)})^{2}]_{n_{1}}$ is uniformly integrable. Both these conditions will be
satisfied if we are able to verify the conditions of Theorem \ref{PU05}, in
the Appendix, for the sequence $(a_{1}Y_{n,1}^{(\ell)}+...+a_{N}Y_{n,N}%
^{(\ell)})_{n}$, where $a_{1},...,a_{N}$ are arbitrary, fixed real numbers. We
have to show that, for $\ell$ fixed $\ $%
\begin{equation}
\sum\nolimits_{k\geq1}\frac{1}{k^{3/2}}||\sum_{j=1}^{k}E(a_{1}Y_{j,1}^{(\ell
)}+...+a_{N}Y_{j,N}^{(\ell)}|\mathcal{F}_{1,N}^{(\ell)})||<\infty.
\label{cramer}%
\end{equation}
By the triangle inequality it is enough to treat each sum separately and to
show that for all $1\leq v\leq N$ we have%
\[
\sum\nolimits_{k\geq1}\frac{1}{k^{3/2}}||\sum_{j=1}^{k}E(Y_{j,v}^{(\ell
)}|\mathcal{F}_{1,N}^{(\ell)})||<\infty.
\]
By (\ref{pcf}) we have that $E(Y_{j,v}^{(\ell)}|\mathcal{F}_{1,N}^{(\ell
)})=E(Y_{j,v}^{(\ell)}|\mathcal{F}_{1,v}^{(\ell)}).$ Therefore, by
stationarity, the latter condition is satisfied if we can prove that%
\[
\sum\nolimits_{k\geq1}\frac{1}{k^{3/2}}||\sum_{j=1}^{k}E(Y_{j,1}^{(\ell
)}|\mathcal{F}_{1,1}^{(\ell)})||<\infty.
\]
Now, by using once again (\ref{pcf}), we deduce%
\begin{align*}
E(Y_{j,1}^{(\ell)}|\mathcal{F}_{1,1}^{(\ell)})  &  =E(X_{j,1}^{(\ell
)}-E(X_{j,1}^{(\ell)}|\mathcal{F}_{j,0}^{(\ell)})|\mathcal{F}_{1,1}^{(\ell
)})=\\
&  E(X_{j,1}^{(\ell)}|\mathcal{F}_{1,1}^{(\ell)})-E(X_{j,1}^{(\ell
)}|\mathcal{F}_{1,0}^{(\ell)}).
\end{align*}
So, by the triangle inequality and the monotonicity of the $L^{2}-$norm of the
conditional expectation with respect to increasing random fields, we obtain%
\[
||\sum_{j=1}^{k}E(Y_{j,1}^{(\ell)}|\mathcal{F}_{1,1}^{(\ell)})||\leq
2||\sum_{j=1}^{k}E(X_{j,1}^{(\ell)}|\mathcal{F}_{1,1}^{(\ell)})||=2\frac
{1}{\ell^{1/2}}||E(S_{k,\ell}|\mathcal{F}_{1,\ell})||.
\]
Furthermore, since the filtration is commuting, by the triangle inequality we
obtain%
\[
||E(S_{k,\ell}|\mathcal{F}_{1,\ell})||=||\sum_{u=1}^{k}\sum_{v=1}^{\ell
}E(X_{u,v}^{\ }|\mathcal{F}_{1,v})||\leq\ell||\sum_{u=1}^{k}E(X_{u,1}%
^{\ }|\mathcal{F}_{1,1})||.
\]
By taking into account condition (\ref{MW}), it follows that we have
\[
\sum\nolimits_{k\geq1}\frac{1}{k^{3/2}}||\sum_{j=1}^{k}E(Y_{j,v}^{(\ell
)}|\mathcal{F}_{1,N}^{(\ell)})||\leq2\ell^{1/2}\sum\nolimits_{k\geq1}\frac
{1}{k^{3/2}}||E(S_{k,1}^{\ }|\mathcal{F}_{1,1})||<\infty,
\]
showing that condition (\ref{cramer}) is satisfied, which implies that the
conditions of Corollary \ref{cormain} are satisfied. The conclusion is that
\[
\frac{1}{\sqrt{n_{1}k}}\sum\nolimits_{j=1}^{n_{1}}\sum\nolimits_{i=1}%
^{k}Y_{j,i}^{(\ell)}\Rightarrow N(0,\sigma_{\ell}^{2})\text{ as }\min
(n_{1},k)\rightarrow\infty,
\]
where $\sigma_{\ell}^{2}$ is defined, in accordance with Theorem \ref{PU05},
by
\[
\sigma_{\ell}^{2}=\lim_{n\rightarrow\infty}\frac{1}{n}E\left(  \sum
\nolimits_{j=1}^{n}Y_{j,1}^{(\ell)}\right)  ^{2}.
\]
According to Theorem 3.2 in Billingsley (1999), in order to prove convergence
and to find the limiting distribution of $S_{n_{1},n_{2}}/\sqrt{n_{1}n_{2}}$
we have to show that
\begin{equation}
\lim_{\ell\rightarrow\infty}\limsup_{n_{1},k\rightarrow\infty}||\frac{1}%
{\sqrt{n_{1}n_{2}}}S_{n_{1},n_{2}}-\frac{1}{\sqrt{n_{1}k}}\sum\nolimits_{j=1}%
^{n_{1}}\sum\nolimits_{i=1}^{k}Y_{j,i}^{(\ell)}||=0 \label{show1}%
\end{equation}
and $N(0,\sigma_{\ell}^{2})\Rightarrow N(0,\sigma^{2}),$ which is equivalent
to%
\begin{equation}
\sigma_{\ell}^{2}\rightarrow\sigma^{2}\text{ as }\ell\rightarrow\infty.
\label{show 2}%
\end{equation}
The conclusion will be that $S_{n_{1},n_{2}}/\sqrt{n_{1}n_{2}}\Rightarrow
N(0,\sigma^{2})$ as $\min(n_{1},n_{2})\rightarrow\infty.$

Let us first prove (\ref{show1}). By the triangle inequality we shall
decompose the difference in (\ref{show1})\ into two parts. Relation
(\ref{show1}) will be established if we show both
\begin{equation}
\lim_{\ell\rightarrow\infty}\limsup_{n_{1},k\rightarrow\infty}\frac{1}%
{\sqrt{n_{1}k}}||\sum\nolimits_{j=1}^{n_{1}}\sum\nolimits_{i=1}^{k}%
E(X_{j,i}^{(\ell)}|\mathcal{F}_{j,i-1}^{(\ell)})||=0. \label{limit}%
\end{equation}
and
\begin{equation}
\lim_{n_{1},k\rightarrow\infty}||\frac{1}{\sqrt{n_{1}n_{2}}}S_{n_{1},n_{2}%
}-\frac{1}{\sqrt{n_{1}k\ell}}S_{n_{1},k\ell}||=0. \label{limit2}%
\end{equation}
In order for computing the standard deviation of the double sum involved,
before taking the limit in (\ref{limit}), we shall apply Lemma \ref{var0} and
a multivariate version of Remark \ref{equivalent} in the Appendix. This
expression is dominated by a universal constant times%
\[
\sum\nolimits_{i,j\geq1}\frac{1}{(ij)^{3/2}}||\sum\nolimits_{u=1}^{j}%
\sum\nolimits_{v=1}^{i}E(E(X_{u,v}^{(\ell)}|\mathcal{F}_{u,v-1}^{(\ell
)})|\mathcal{F}_{1,0}^{(\ell)})||.
\]
Now,
\[
\sum\nolimits_{u=1}^{j}\sum\nolimits_{v=1}^{i}E(E(X_{u,v}^{(\ell)}%
|\mathcal{F}_{u,v-1}^{(\ell)})|\mathcal{F}_{1,0}^{(\ell)})=\frac{1}{\ell
^{1/2}}E(S_{j,i\ell}|\mathcal{F}_{1,0}).
\]
So, the quantity in (\ref{limit}) is bounded above by a universal constant
times
\[
\frac{1}{\ell^{1/2}}\sum\nolimits_{i,j\geq1}\frac{1}{(ij)^{3/2}}%
||E(S_{j,i\ell}|\mathcal{F}_{1,0})||,
\]
which converges to $0$ as $\ell\rightarrow\infty$ under our condition
(\ref{MW}), by Lemmas 2.7 and 2.8 in Peligrad and Utev (2005), applied in the
second coordinate.

As far as the limit (\ref{limit2}) is concerned, since by Lemma \ref{var0} and
condition (\ref{MW}) the array $\sum\nolimits_{j=1}^{n_{1}}\sum\nolimits_{i=1}%
^{n_{2}}X_{j,i}/\sqrt{n_{1}\ n_{2}}$ is bounded in $L^{2}$, it is enough to
show that, for $k\ell<n_{2}<(k+1)\ell,$ we have
\[
\lim_{n_{1},n_{2}\rightarrow\infty}||\frac{1}{\sqrt{n_{1}n_{2}}}%
\sum\nolimits_{j=1}^{n_{1}}\sum\nolimits_{i=k\ell+1}^{n_{2}}X_{j,i}||=0.
\]
We just have to note that, again by Lemma \ref{var0}, condition (\ref{MW}) and
stationarity, there is a constant $K$ such that
\[
||\sum\nolimits_{j=1}^{n_{1}}\sum\nolimits_{i=k\ell+1}^{n_{2}}X_{j,i}||\leq
K\sqrt{n_{1}\ell}%
\]
and $\ell/n_{2}\rightarrow0$ as $n_{2}\rightarrow\infty.$

We turn now to prove (\ref{show 2})$.$ By (\ref{show1}) and the orthogonality
of martingale differences,%
\[
\lim_{\ell\rightarrow\infty}\limsup_{n_{1},n_{2}\rightarrow\infty}|\text{
}||\frac{1}{\sqrt{n_{1}n_{2}}}S_{n_{1},n_{2}}||-||\frac{1}{\sqrt{n_{1}}}%
\sum\nolimits_{j=1}^{n_{1}}Y_{j,0}^{(\ell)}||\text{ }|=0.
\]
So%
\[
\lim_{\ell\rightarrow\infty}\limsup_{n_{1},n_{2}\rightarrow\infty}|\text{
}||\frac{1}{\sqrt{n_{1}n_{2}}}S_{n_{1},n_{2}}||-\sigma_{\ell}|=0.
\]
By the triangle inequality, this shows that $\sigma_{\ell}$ is a Cauchy
sequence, therefore convergent to a constant $\sigma$ and also%
\[
\lim_{n_{1},n_{2}\rightarrow\infty}||\frac{1}{\sqrt{n_{1}n_{2}}}S_{n_{1}%
,n_{2}}||=\sigma.
\]
The proof is now complete. $\square$

\bigskip

The extensions to random fields indexed by $Z^{d},$ for $d>2,$ are
straightforward following the same lines of proofs as for a two-indexed random
field. We shall point out the differences. To extend Lemma\textbf{ }%
\ref{var0}, we first apply a result of Peligrad and Utev (2005) (see Theorem
\ref{PU05} in the Appendix) to the stationary sequence $Y_{\mathbf{j}}%
(m)=\sum\nolimits_{i=1}^{m}X_{\mathbf{j},i}$ with $\mathbf{j\in}Z^{d-1}$ and
then we apply induction.

In order to prove Theorem \ref{general}, we partition the variables according
to the last index. Let $\ell$ be a fixed positive integer, denote
$k=[n_{d}/\ell]$ and define
\[
X_{\mathbf{j},i}^{(\ell)}=\frac{1}{\ell^{1/2}}\sum_{u=(i-1)\ell+1}^{i\ell
}X_{\mathbf{j},u}\;,i\geq1.
\]
Then, for each $\mathbf{j}$ we construct the stationary sequence of martingale
differences $(Y_{\mathbf{j},i}^{(\ell)})_{i\in Z}$ defined by $Y_{\mathbf{j}%
,i}^{(\ell)}=X_{\mathbf{j},i}^{(\ell)}-E(X_{\mathbf{j},i}^{(\ell)}%
|\mathcal{F}_{\mathbf{j},i-1}^{(\ell)})$ and
\[
D_{\mathbf{n}^{\prime},i}^{(\ell)}=\frac{1}{\sqrt{|\mathbf{n}^{\prime}|}}%
\sum\nolimits_{\mathbf{j}=1}^{\mathbf{n}^{\prime}}Y_{\mathbf{j},i}^{(\ell)}.
\]
For showing that $(D_{\mathbf{n}^{\prime},1}^{(\ell)},...,D_{\mathbf{n}%
^{\prime},N}^{(\ell)})\Rightarrow(L_{1},...,L_{N}),$ we apply the induction
hypothesis. $\ \square$

\bigskip

\textbf{Proof of Example \ref{exinear}}.

Let us note first that the variables are square integrable and well defined.
Note that%
\[
E(S_{\mathbf{u}}|\mathcal{F}_{\mathbf{0}})=\sum_{\mathbf{1}\leq\mathbf{k\leq
u}}\sum_{\mathbf{j}\leq\mathbf{0}}a_{\mathbf{k}-\mathbf{j}}\xi_{\mathbf{j}}%
\]
and therefore%
\[
E(E^{2}(S_{\mathbf{u}}|\mathcal{F}_{\mathbf{0}}))=\sum_{\mathbf{i}%
\geq\mathbf{0}}(\sum_{\mathbf{1}\leq\mathbf{k\leq u}}a_{\mathbf{k}+\mathbf{i}%
})^{2}E(\xi_{\mathbf{1}}^{2}).
\]
The result follows by applying Theorem \ref{general} (see Remark
\ref{equivalent} and consider a multivariate analog of it). $\square$

\bigskip

\textbf{Proof of Example \ref{Volterra}}.

Note that%
\begin{gather*}
E(S_{\mathbf{j}}|\mathcal{F}_{\mathbf{0}})=\sum\limits_{\mathbf{k}=\mathbf{1}%
}^{\mathbf{j}}\sum_{(\mathbf{u},\mathbf{v)}\geq(\mathbf{k},\mathbf{k}%
)}a_{\mathbf{u},\mathbf{v}}\xi_{\mathbf{k-u}}\xi_{\mathbf{k-v}}\\
=\sum_{(\mathbf{u},\mathbf{v)}\geq(\mathbf{0},\mathbf{0})}\sum
\limits_{\mathbf{k}=\mathbf{1}}^{\mathbf{j}}a_{\mathbf{k+u},\mathbf{k+v}}%
\xi_{-\mathbf{u}}\xi_{-\mathbf{v}}=\sum_{(\mathbf{u},\mathbf{v)}%
\geq(\mathbf{0},\mathbf{0})}c_{\mathbf{u},\mathbf{v}}(\mathbf{j}%
)\xi_{-\mathbf{u}}\xi_{-\mathbf{v}}.
\end{gather*}
Since by our conditions $c_{\mathbf{u},\mathbf{u}}=0$ we obtain
\[
E(E^{2}(S_{\mathbf{j}}|\mathcal{F}_{\mathbf{0}}))=\sum\limits_{\mathbf{u}%
\geq\mathbf{0,v}\geq\mathbf{0,u}\neq\mathbf{v}}(c_{\mathbf{u},\mathbf{v}}%
^{2}(\mathbf{j})+c_{\mathbf{u},\mathbf{v}}(\mathbf{j})c_{\mathbf{v}%
,\mathbf{u}}(\mathbf{j}))E(\xi_{\mathbf{u}}\xi_{\mathbf{v}})^{2}.
\]
$\square$

\section{Appendix.}

For convenience we mention a classical result of McLeish which can be found on
pp. 237-238 G\"{a}nssler and H\"{a}usler (1979).

\begin{theorem}
\label{martCLT}Assume $(D_{n,i})_{1\leq i\leq n}$ is an array of square
integrable martingale differences adapted to an array $(\mathcal{F}%
_{n,i})_{1\leq i\leq n}$ of nested sigma fields. Suppose that
\[
\max_{1\leq j\leq n}|D_{n,j}|\rightarrow^{L_{2}}0\text{ as }n\rightarrow
\infty.
\]
and%
\[
\sum_{j=1}^{n}D_{n,j}^{2}\rightarrow^{P}c^{2}\ \text{ }\ \text{as
}n\rightarrow\infty\text{.}%
\]
Then $\sum_{j=1}^{n}D_{n,j}$ converges in distribution to $N(0,c^{2}).$
\end{theorem}

The following is a Corollary of Theorem 1.1 in Peligrad and Utev (2005). This
central limit theorem was obtained by Maxwell and Woodroofe (2000).

\begin{theorem}
\label{PU05}Assume that $(X_{i})_{i\in Z}$ is a stationary sequence adapted to
a stationary filtration $(\mathcal{F}_{i})_{i\in Z}$. Then there is a
universal constant $C_{1}$ such that
\[
||S_{n}||\leq C_{1}n^{1/2}\sum\nolimits_{k=1}^{\infty}\frac{1}{k^{3/2}%
}||E(S_{k}|\mathcal{F}_{1})||.
\]
If%
\[
\sum\nolimits_{k=1}^{\infty}\frac{1}{k^{3/2}}||E(S_{k}|\mathcal{F}%
_{1})||<\infty,
\]
then $(S_{n}^{2}/n)_{n}$ is uniformly integrable and and there is a positive
constant $c$ such that
\[
\frac{1}{n}E\left(  S_{n}\right)  ^{2}\rightarrow c^{2}\text{ as }%
n\rightarrow\infty.
\]
If in addition the sequence is ergodic we have
\[
\frac{1}{\sqrt{n}}S_{n}\Rightarrow cN(0,1)\text{ as }n\rightarrow\infty.
\]

\end{theorem}

\begin{remark}
\label{equivalent}Note that we have the following equivalence:%
\[
\sum\nolimits_{k=1}^{\infty}\frac{1}{k^{3/2}}||E(S_{k}|\mathcal{F}%
_{1})||<\infty\text{ if and only if }\sum\nolimits_{k=1}^{\infty}\frac
{1}{k^{3/2}}||E(S_{k}|\mathcal{F}_{0})||<\infty.
\]

\end{remark}

\begin{remark}
\label{remPU05}The condition (\ref{MW1}) is implied by%
\[
\sum\nolimits_{k=1}^{\infty}\frac{1}{k^{1/2}}||E(X_{k}|\mathcal{F}%
_{1})||<\infty.
\]

\end{remark}

\begin{lemma}
\label{commuting}Assume that $X,Y,Z$ are integrable random variables such that
$(X,Y)$ and $Z$ are independent. Assume that $\ g(X,Y)$ is integrable. Then%
\[
E(g(X,Y)|\sigma(Y,Z))=E(g(X,Y)|Y)\text{ a.s.}%
\]
and%
\[
E(g(Z,Y)|\sigma(X,Y))=E(g(Z,Y)|Y)\text{ a.s.}%
\]

\end{lemma}

Proof. Since $(X,Y)$ and $Z$ are independent, it is easy to see that $X$ and
$Z$ are conditionally independent given $Y$. The result follows from this
observation by Problem 34.11 in Billingsley (1995). $\square$

\bigskip

\textbf{Acknowledgements.} This research was supported in part by the
NSF\ grant DMS-1512936 and the Taft Research Center at the University of
Cincinnati. The authors would like to thank Yizao Wang for pointing out that
filtrations generated by random fields with independent columns generate
commuting filtrations. We would also like to thank anonymous referees for
carefully reading the manuscript and their numerous suggestions, which
improved the presentation of the paper.


\begin{thebibliography}{99}                                                                                               %


\bibitem {B}Billingsley, P. (1995). Probability and measures. (3rd ed.). Wiley
Series in Probability and Statistics, New York

\bibitem {B2}Billingsley, P. (1999). Convergence of probability measures. (2nd
ed.). Wiley Series in Probability and Statistics, New York.

\bibitem {Br}Bradley, R. and Tone, C. (2017). A Central Limit Theorem for
Non-Stationary Strongly Mixing Random Fields. Journal of Theoretical
Probability 30 655-674.

\bibitem {CM}Cuny, C. and Merlev\`{e}de, F. (2014). On martingale
approximations and the quenched weak invariance principle. Ann. Probab. 42 760-793.

\bibitem {CDV}Cuny, C. Dedecker J. and Voln\'{y}, D. (2015). A functional
central limit theorem for fields of commuting transformations via martingale
approximation, \textit{Zapiski Nauchnyh Seminarov} POMI 441.C. Part 22 239-263
and Journal of Mathematical Sciences 2016, 219 765--781.

\bibitem {D}Dedecker, J. (1998). A central limit theorem for stationary random
fields. Probability Theory and Related Fields. 110 397-426. \ 

\bibitem {EM}El Machkouri, M., Voln\'{y}, D. and Wu, W.B. (2013). A central
limit theorem for stationary random fields. Stochastic Process. Appl. 123 1-14.

\bibitem {GH}G\"{a}nssler, P. and H\"{a}usler, E. (1979). Remarks on the
Functional Central Limit Theorem for Martingales. Z.
Wahrscheinlichkeitstheorie verw. Gebiete 50 237-243.

\bibitem {G1}Gordin M. I. (1969). On the central limit theorem for stationary
processes. Dokl. Akad. Nauk SSSR, 188 739--741. and Soviet Math. Dokl. 10 1174--1176.

\bibitem {G}Gordin, M. I. (2009). Martingale co-boundary representation for a
class of stationary random fields, Zap. Nauchn. Sem. S.-Peterburg. Otdel. Mat.
Inst. Steklov. (POMI) 364, Veroyatnostn i Statistika. 14.2, 88-108, 236; and
J. Math. Sci. 163 (2009) 4, 363-374.

\bibitem {mw}Maxwell, M. and Woodroofe, M. (2000). Central limit theorems for
additive functionals of Markov chains. Ann. Probab. 28, 713--724.

\bibitem {Pu05}Peligrad, M. and Utev, S. (2005). A new maximal inequality and
invariance principle for stationary sequences. Ann. Probab. 33, 798-815.

\bibitem {R}Rosenblatt, M. (1972). Central limit theorem for stationary
processes, Berkeley Symp. on Math. Statist. and Prob. Proc. Sixth Berkeley
Symp. on Math. Statist. and Prob., Vol. 2 (Univ. of Calif. Press), 551-561.

\bibitem {V}Voln\'{y}, D. (2015). A central limit theorem for fields of
martingale differences, C. R. Math. Acad. Sci. Paris 353, 1159-1163.

\bibitem {VW}Voln\'{y}, D. and Wang, Y. (2014). An invariance principle for
stationary random fields under Hannan's condition. Stochastic Proc. Appl. 124 4012-4029.

\bibitem {WW}Wang Y. and Woodroofe, M. (2013). A new criteria for the
invariance principle for stationary random fields. Statistica Sinica 23 1673-1696.
\end{thebibliography}
\end{document}